\documentclass[12pt]{article}

\usepackage{cmap}
\usepackage[T2A]{fontenc}
\usepackage[utf8]{inputenc} 
\usepackage[english]{babel}
\usepackage{amsmath, amssymb, latexsym}
\usepackage{tikz}
\usepackage{ifthen}
\usepackage{hyperref}
\usepackage{longtable}
\usepackage{setspace}
\title{Extended perfect codes}
\date{}

\makeatletter
\def\@seccntformat#1{\csname the#1\endcsname.\ } 
\makeatother

\newif\ifNoRemark
\def\addtheorem#1#2#3#4{
\ifthenelse{\equal{#2}{}}{}%
{\ifthenelse{\expandafter\isundefined\csname the#2\endcsname}{\newcounter{#2}}{}}
\newenvironment{#1}[1][\global\NoRemarktrue]
{\par\addvspace{2mm plus 0.5mm minus 0.2mm}\noindent 
{\bf #3}\ifthenelse{\equal{#2}{}}{}%
{\refstepcounter{#2}{\bf ~\csname the#2\endcsname}}%
{\bf \vphantom{##1}\ifNoRemark.\ \else\ (##1).\fi}\begingroup #4}%
{\endgroup\par\addvspace{1mm plus 0.5mm minus 0.2mm}\global\NoRemarkfalse}
\expandafter\newcommand\csname b#1\endcsname{\begin{#1}}
\expandafter\newcommand\csname e#1\endcsname{\end{#1}}
}

\addtheorem{theoreman}{thrm}{Theorem}{\sl}
\addtheorem{lemman}{lmm}{Lemma}{\sl}
\addtheorem{predln}{prpstn}{Proposition}{\sl}
\addtheorem{coroll}{crllr}{Corollary}{\sl}
\addtheorem{remarkn}{rmrk}{Remark}{}
\addtheorem{definitionn}{def}{Definition}{}

 \newenvironment{proof}[1][\hspace{-1.0ex}]%
  {\par\addvspace{1mm}{\sc Proof\hspace{1.0ex}{#1}.} }%
  {\quad$\blacktriangle$\par\addvspace{1mm}}
	
\title{On the non-existence of extended perfect codes and some perfect colorings%
\thanks{This work was funded by the Russian Science Foundation (Grant 18-11-00136).} 
}
\author{Evgeny Bespalov
\thanks{Sobolev Institute of Mathematics, Novosibirsk, Russia. E-mail: bespalovpes@mail.ru}
}

\begin{document}
\maketitle
\def\VV{{\scriptscriptstyle\mathrm{V}}}

\begin{abstract}
    In this paper we obtain the necessary condition for the existence of perfect $k$-colorings (equitable $k$-partitions) 
    in Hamming graphs $H(n,q)$, where $q=2,3,4$ and Doob graphs $D(m,n)$.  
    As an application, we prove the non-existence of extended perfect codes in $H(n,q)$, 
    where $q=3,4$, $n>q+2$, and in $D(m,n)$, where $2m+n>6$.
\end{abstract}

\section{Introduction}

A $k$-coloring of a graph $G=(V,E)$ 
is a surjective function from the vertex set $V$ into a color set of cardinality $k$, usually denoted by $\{0,1,\ldots,k-1\}$. 
This coloring is called perfect if for any $i,j$ the number of vertices of color $j$ in the neighbourhood of vertex $x$ of color $i$ depends only on $i$ and $j$, 
but not on the choice of $x$.
An equivalent concept is an equitable $k$-partition, which is a partition of the vertex set $V$ into cells $V_0,\ldots,V_{k-1}$, 
where these cells are the preimages of the colors of some perfect $k$-coloring. 
Also the perfect colorings are the particular cases of the perfect structures, see e.g. \cite{Tar:perfstruct}.
In this paper, we consider perfect colorings in Hamming graphs $H(n,q)$ 
(mainly focusing on the case $q=2,3,4$) and Doob graphs $D(m,n)$.
Remind that the Hamming graph $H(n,q)$ is the direct product of $n$ copies of the complete graph $K_q$ on $q$ vertices, and 
the Doob graph $D(m,n)$, where $m>0$, is the direct product of $m$ copies of the Shrikhande graph and $n$ copies of $K_4$.
These graphs are distance-regular; moreover, the Doob graph $D(m,n)$ has the same intersection array as $H(2m+n,4)$.
Many combinatorial objects can be defined as perfect colorings with corresponding parameters,
for example, MDS codes with distance $2$; 
latin squares and latin hypercubes; 
unbalanced boolean functions attending the correlation-immunity bound \cite{FDF:CorrImmBound};
orthogonal arrays attaining the Bierbrauer--Friedman bound \cite{Bierbrauer:95,Friedman:92};
boolean-value functions on Hamming graphs and orthogonal arrays that attach some other bounds \cite{Pot:2012:color,Pot:2010:correng,Kro:OA1536};
some binary codes attending the linear-programming bound that are cells of equitable partitions into 4, 5, or 6 cells \cite{Kro:2m-3,Kro:2m-4}.

One important class of objects that corresponding to perfect colorings is the $1$-perfect codes.
It is generally known \cite[Ch.~6, Th.~37]{MWS} that if $q=p^m$ is a prime power,
then there is a $1$-perfect code in $H(n,q)$ if and only if $n=(q^l-1)/(q-1)$ for some positive integer $l$.
In the case when $q$ is not prime power, there is \emph{very little known} about the existence of $1$-perfect codes. 
It is known that there are no $1$-perfect codes in $H(7,6)$ \cite[Theorem~6]{GolombPosner64} (since there are no pair of orthogonal latin squares of order $6$). 
Heden and Roos obtained the necessary condition \cite{HedRoos:2011} on the non-existence of some $1$-perfect codes, 
which in particular implies the non-existence of $1$-perfect codes in $H(19,6)$. 
Also we mention result of Lenstra \cite{Lenstra72}, 
which generalized Lloyd's condition (see \cite{Lloyd,MWS}) for a non-prime power $q$.  
This result implies that if there is a $1$-perfect code in $H(n,q)$, then $n=kq+1$.
Krotov \cite{Kro:pfdoob} completely solved the problem of the existence of $1$-perfect codes in Doob graphs. 
Namely, he proved that there is a $1$-perfect code in $D(m,n)$ if and only if $2m+n=(4^l-1)/3$ for some positive integer $l$. 
Note that the existence of a $1$-perfect code in $D(m,n)$ not always implies 
the existence of linear or additive $1$-perfect codes in this graph 
(the set of admissible parameters of unrestricted $1$-perfect codes in Doob graphs is essentially wider than that of linear \cite{Kro:perfect-doob} or additive \cite{SHK:addperfdoob} $1$-perfect codes).

Another important class of codes corresponding to perfect colorings is the completely regular codes.
A code $C$ is completely regular if the distance coloring with respect to $C$ 
(a vertex $v$ has the color that is equal to the distance from $v$ to $C$) is perfect. 
These codes originally were defined by Delsarte \cite{Delsarte:1973}, but here we use the different equivalent definition from \cite{Neu:crg}.
For more information about completely regular codes and problem of its existence, 
we refer to the survey \cite{BRZ:crg}, papers \cite{FDF:PerfCol, BKMTV:perfcolinham} (for codes with covering radius $\rho=1$), and
the small-value tables of parameters \cite{KKM:smallval}.

In the current paper, we \emph{stay} on the class of completely regular codes that correspond to extended $1$-perfect codes. 
An extended $1$-perfect code is a code with code distance $4$ 
obtained by appending an additional symbol (this operation we call an extension) to the codewords of some $1$-perfect code  
(the rigorous definition will be given in the next section).
It is known that there is an extended $1$-perfect code in $H(2^m,2)$ and in $H(2^m+2,2^m)$ 
for any positive integer $m$ (see~\cite{MWS, BRZ:crg}).
It was mentioned in \cite[Section~4]{AhlAydKha} that a result from \cite{Hill:Caps}
implies the non-existence of extended $1$-perfect codes in $H(n,q)$ obtained from the Hamming codes, 
except the case when $(n,q)=(2^m,2)$ or $(n,q)=(2^m+2,2^m)$. 
An extended $1$-perfect code in $H(q+2,q)$ (or in $D(m,n)$, where $2m+n=6$) is also an MDS code with distance $4$.
There is a characterization of all extended $1$-perfect codes in $H(6,4)$ \cite{Alderson:MDS4} and in Doob graphs $D(m,n)$ \cite{BesKro:mdsdoob}, where $2m+n=6$, $m>0$. 
Ball showed \cite{Ball:2012:1} that if $q$ is odd prime, then there are no linear extended $1$-prefect codes in $H(q+2,q)$. 
The non-existence of extended $1$-perfect codes in $H(7,5)$ and $H(9,7)$ was proved in \cite{KKO:smallMDS}. 
In \cite{KokOst:further} it was shown that any extended $1$-perfect code in $H(10,8)$ is equivalent to a linear code.
The non-existence of extended $1$-perfect codes in $H(14,3)$ follows from the bound established in \cite{GST:newupbounds}.
For completeness, note that formally codes consisting from one vertex in $H(2,q)$ are also extended $1$-perfect codes.
Such codes are called trivial.

In this paper, we obtain a necessary condition for the existence of perfect colorings 
in Hamming graphs $H(n,q)$, where $q=2,3,4$, and Doob graphs.  
We apply it to extended $1$-perfect codes and prove that there are no such codes 
in $H(n,q)$, $q=3,4$, $n>q+2$, and in $D(m,n)$, $2m+n>6$.  
This completes the characterization of such codes in these graphs (see Theorem~\ref{t:parameters}).
In addition, we prove that extended $1$-perfect codes can exist in $H(n,q)$ only if $n$ is even, 
which particularly implies the non-existence of some MDS codes with distance $4$.  
We hope that this method can be applied for a proof of the non-existence of some other perfect $k$-colorings  
(but for perfect $2$-colorings it does not add something new to results from \cite{BKMTV:perfcolinham}).

The paper is organized as follows.  
In Section~\ref{s:prelim}, we give main definitions and simple observations.
In Section~\ref{s:necessary}, we obtain a necessary condition (Theorem~\ref{t:nessesary}) for the existence 
of perfect colorings in Doob graphs and Hamming graphs $H(n,q)$, where $q=2,3,4$.
In Section~\ref{s:extarecrg}, we prove that any extended $1$-perfect code in $H(n,q)$ 
is a completely regular code with intersection array $(n(q-1),(n-1)(q-1);1,n)$, and vise versa;
similar results are shown for Doob graphs.
This allows us to apply Theorem~\ref{t:nessesary} to prove the non-existence of some extended $1$-perfect codes in Section~\ref{s:nonexist}.
Finally, we describe all parameters for which there is an extended $1$-perfect code in $D(m,n)$ and $H(n,3)$ in Theorem~\ref{t:parameters}.

\section{Preliminaries}\label{s:prelim}

Given a graph $G$, we denote by $\VV{G}$ its vertex set.
A surjective function $f: \VV{G} \to \{0,1,\ldots,k-1\}$ on the vertex set of $G$
is called a \emph{$k$-coloring} of a graph $G$ in the colors $0,1,\ldots,k-1$. 
If for all $i,j$ every vertex $x$ of color $i$ has exactly $s_{i,j}$ neighbours of color $j$,
where $s_{i,j}$ does not depend on the choice of $x$, then the coloring $f$ is called a \emph{perfect $k$-coloring} with \emph{quotient matrix} $S=(s_{i,j})$. 

Let $G$ be a connected graph. A \emph{code} $C$ in $G$ is an arbitrary nonempty subset of $\VV{G}$.
The \emph{distance} $d(x,y)$ between two vertices $x$ and $y$ is the length of the shortest path between $x$ and $y$. 
The \emph{code distance} $d$ of a code $C$ is the minimum distance between two different vertices of $C$. 
The distance $d(A,B)$ between two sets of vertices $A$ and $B$ 
equals $\min\{d(x,y):x \in A, y \in B \}$. 
The \emph{covering radius} of a code $C$ is $\rho=\max\limits_{v \in \VV{G}}\{d(\{v\}, C)\}$.
Let $C$ be a code in a graph $G$. 
The \emph{distance coloring} with respect to $C$ is the coloring $f$ defined in the following way: $f(x)$ is equal to the distance between $\{x\}$ and $C$. If $f$ is a perfect coloring with quotient matrix $S$, then $C$ is called a \emph{completely regular code} with quotient matrix $S$.  
In this case, the matrix $S$ is tridiagonal.  
A connected regular graph $G$ is called \emph{distance-regular} if for any vertex $x$ of $G$ the set $\{x\}$ 
is a completely regular code with quotient matrix $S$ that does not depend on the choice of $x$. The sequence
$(b_0,\ldots,b_{\rho-1};c_1,\ldots,c_{\rho})=(s_{0,1},\ldots,s_{\rho-1,\rho};s_{1,0},\ldots,s_{\rho,\rho-1})$ is called the \emph{intersection array}.

The \emph{Shrikhande graph} $Sh$ is a Cayley graph with the vertex set $\mathbb Z^2_4$, 
where two vertices $x$ and $y$ are adjacent if and only if their difference $(x-y)$ belongs to the connecting set $\{01,10,03,30,11,33\}$.
The complete graph $K_q$ on $q$ vertices can be represented as a Cayley graph, where the vertex set is $\mathbb Z_q$,
and two vertices $x$ and $y$ are adjacent if and only if their difference $(x-y)$ belongs to the connecting set $\{1,2,\ldots,q-1\}$.
The Hamming graph $H(n,q)$ is the direct product $K^n_q$ of $n$ copies of $K_q$. 
The vertex set of $H(n,q)$ can be represented as $\mathbb Z^n_q=\{(x_1,\ldots,x_n):x_i \in \mathbb Z_q\}$.
Denote by $D(m,n)=Sh^m \times K^n_4$ the direct product of $m$ copies of the Shrikhande graph $Sh$ and $n$ copies of the complete graph $K_4$.
If $m>0$, then this graph is called \emph{Doob graph}.
The vertex set of $D(m,n)$ can be represented as $(\mathbb Z^2_4)^m \times \mathbb Z^n_4=\{(x_1,\ldots,x_m;y_1,\ldots,y_n): x_i \in \mathbb Z^2_4, y_j \in \mathbb Z_4\}$.
The Hamming graph $H(n,q)$ is distance-regular with intersection array $(n(q-1),(n-1)(q-1),\ldots,q-1;1,2,\ldots,n)$.
The Doob graph $D(m,n)$ is distance-regular with the same intersection array as $H(2m+n,4)$.


For a vertex $v=(x_1,\ldots,x_{n-1})$ of $H(n-1,q)$ and $a \in \mathbb Z_q $, denote by $v^a_i=(x_1,\ldots,x_{i-1},a,x_{i},\ldots,x_{n-1})$ the vertex of $H(n,q)$.
Analogously, for a vertex $v=(x_1,\ldots,x_m;y_1,\ldots,y_{n-1})$ of $D(m,n-1)$ and $a \in \mathbb Z_4$, 
denote by \linebreak $v^a_{;i}=(x_1,\ldots,x_m;y_1,\ldots,y_{i-1},a,y_{i},\ldots,y_{n-1})$ the vertex of $D(m,n)$. 
The \emph{projection} (also known as puncturing) $C_i$ of a code $C$ in $H(n,q)$ is the code in $H(n-1,q)$ defined as follows: 
$$C_i=\{v \in \VV{H(n-1,q)}: v^a_i \in C \text{ for some } a \in \mathbb Z_q\}.$$
Similarly, the projection $C_{;i}$ of a code $C$ in the  Doob graph $D(m,n)$, $n>0$, is the code in $D(m,n-1)$ defined as follows 
$$C_{;i}=\{v \in \VV{D(m,n-1)}: v^a_{;i} \in C \text{ for some } a \in \mathbb Z_4\}.$$

By $B_e(x)=\{y:d(x,y) \le e\}$ denote the radius $e$-ball with center $x$. 
A code $C$ in a graph $G$ is called \emph{$e$-perfect} if $|C \cap B_e(x)|=1$ for any $x \in  \VV{G}$. 
In equivalent definition, an $e$-perfect code is a code with code distance $d=2e+1$, 
whose cardinality achieves the sphere-packing bound.
It is known that if $q=p^m$ is a prime power, then a $1$-perfect code in $H(n,q)$ exists 
if and only if $n=(q^l-1)/(q-1)$  for some positive integer $l$ \cite{MWS}. 
The cardinality of this code is equal to $q^{n-l}$.
It is also known that a $1$-perfect code in $D(m,n)$  exists if and only if $2m+n=(4^l-1)/3$ for some positive integer $l$ \cite{Kro:pfdoob}.
The cardinality of this code is equal to $4^{2m+n-l}$. 

A code $C$ in $H(n,q)$ is called an \emph{extended $1$-perfect code} if its code distance is equal to $4$
and the projection $C_i$ in some position $i$ is a $1$-perfect code.
If a code $C$ has distance $d>1$, then its projection has distance at least $d-1$ and the same cardinality. 
Therefore, if $C$ is an extended $1$-perfect code, then  the projection $C_i$ is a $1$-perfect code for any $i=1,\ldots,n$.
So, if $q=p^m$ is a prime power, then an extended $1$-perfect code in $H(n,q)$ can exist 
only for $n=(q^l+q-2)/(q-1)$, $l \in \mathbb N$. The cardinality of such code equals $q^{n-l-1}$.
Similarly, a code $C$ in $D(m,n)$ is called an \emph{extended $1$-perfect code} if its code distance equals $4$ and the projection $C_{;i}$ for some position $i$ is a $1$-perfect code. 
So an extended $1$-perfect code in $D(m,n)$ can exist only if $2m+n=(4^l+2)/3$, $l \in \mathbb N$.
The cardinality of such code equals $4^{2m+n-l-1}$.
If $n=0$, then a code $C$ in $D(m,0)$ is called an \emph{extended $1$-perfect code} 
if it has the same parameters as an extended $1$-perfect code in Doob graph of the same diameter,
i.e. $2m=(4^l+2)/3$, the code distance is equal to $4$, $|C|=4^{2m-l-1}$.

\section{A necessary condition for the existence of perfect colorings}\label{s:necessary}

Given a graph $G$, let us consider the set of complex-valued functions $f: \VV{G} \to \mathbb C$ on the vertex set. 
These functions form a vector space $U(G)$ with the inner product 
$(f,g)=\sum\limits_{x \in \VV{G}} f(x)\overline{g(x)}$. 
A function $f:\VV{G} \to \mathbb C$ is called an \emph{eigenfunction} of $G$ if $Mf=\lambda f$, $f \not\equiv 0$,
where $M$ is the adjacency matrix of $G$, for some $\lambda$, which is called an \emph{eigenvalue} of $G$.  
Denote by $U_{\lambda}=\{f:Mf=\lambda f\}$ the eigensubspace corresponding to $\lambda$.

Let $G$ be a Hamming graph $H(n,q)$ or a Doob graph $D(m,n)$.
Then it is convenient to use the characters to form a basis of each eigensubspace.
Let $\xi$ be the $q$-th root of unity, namely 
$\xi=e^{\frac{2\pi \sqrt{-1}}{q}}$. 
If $G$ is $H(n,q)$, then for an arbitrary $z \in \mathbb Z^n_q$ define the function $\varphi_z(t)=\frac{\xi^{\langle z,t \rangle}}{q^{n/2}}$, 
where $\langle v,u \rangle =v_1u_1+\ldots+v_nu_n \mod q$.
If $G$ is $D(m,n)$, then for an arbitrary $z \in  (\mathbb Z^2_4)^m \times \mathbb Z^n_4$ define the function $\varphi_z(t)=\frac{\xi^{\langle z,t \rangle}}{4^{(2m+n)/2}}$, 
where $\langle x,v \rangle =(x_{1}v_{1}+y_{1}u_{1})+\ldots+(x_{m}v_{m}+y_{m}u_{m})+r_1s_1+\ldots+r_ns_n \mod 4$; $x=([x_{1},y_{1}],\ldots,[x_{m},y_{m}];r_1,\ldots,r_n)$ and $v=([v_{1},u_{1}],\ldots,[v_{m},u_{m}];s_1,\ldots,s_n)$ are vertices in $D(m,n)$ (we denote by $[a,b]$ an element of $\mathbb Z^2_4$).
It is known that the functions $\varphi_z$, where $z \in \mathbb \VV{G}$, are eigenfunctions of $G$ and these functions
form an orthonormal basis of the vector space $U(G)$. 

\begin{lemman}\label{l:mffs}
Let $f$ be a perfect $k$-coloring of a graph $G$ with quotient matrix $S$. 
Let $f_j=\chi_{f^{-1}(j)}$ be the characteristic function of the set of vertices of color $j$. 
Then for any $t \in \mathbb N$ $$(M^tf_j,f_j)=s^t_{j,j} \cdot |f^{-1}(j)|,$$ where $M$ is the adjacency matrix of $G$, 
and $s^t_{j,j}$ is the $(j,j)$-th element of the matrix $S^t$.
\end{lemman}
\begin{proof}
Let $F=(f_0,\ldots,f_{k-1})$ be the $|\VV{G}| \times k$ matrix, where the $i$-th column $f_i=\chi_{f^{-1}(i)}$ is the characteristic function of the set of vertices of color $i$. 
It is known that $MF=FS$ (see for example \cite[Section~5.2]{Godsil93}), and consequently $M^tF=FS^t$ for any $t$. 
Hence, $(M^tf_j)(x)=(M^tF)_{x,j}=(FS^t)_{x,j}=s^t_{f(x),j}$ for any vertex $x \in \VV{G}$. Since $f_j(x)=0$ if $f(x) \ne j$, we have
$(M^tf_j,f_j)=s^t_{j,j} \cdot |f^{-1}(j)|$.
\end{proof}

\begin{lemman}\cite[Section.~5.2]{Godsil93}.\label{l:eigenvalues}
Let $f$ be a perfect coloring of a graph $G$ with quotient matrix $S$. 
If $\lambda$ is an eigenvalue of $S$, then $\lambda$ is an eigenvalue of $G$.
\end{lemman}



\begin{theoreman}\label{t:nessesary}
Let $G$ be the Hamming graph $H(n,q)$, where $q \in \{2,3,4\}$, or 
the Doob graph $D(m,n)$. Let $f$ be a perfect $k$-coloring of $G$ with quotient matrix $S$ that has eigenvalues 
$\lambda_0 > \lambda_1 > \ldots > \lambda_{l}$. 
Let $i$ be a color of $f$, and let $s^t_{i,i}$ be the $(i,i)$-th element of $S^t$, $t=1,\ldots,l-1$.
  

Then the linear system of equations
\[\displaystyle{\begin{pmatrix}
1 & 1 & \ldots & 1 \\
\lambda_1  & \lambda_2 & \ldots & \lambda_l \\ 
\lambda^2_1  & \lambda^2_2 & \ldots & \lambda^2_l \\
\ldots & \ldots & \ldots & \ldots \\ 
\lambda^{l-1}_1  & \lambda^{l-1}_2 & \ldots & \lambda^{l-1}_l
\end{pmatrix}
\begin{pmatrix}
x_1       \\
x_2  \\ 
x_2  \\
\ldots        \\ 
x_{l}
\end{pmatrix}
= |f^{-1}(i)| 
\begin{pmatrix}
1\\ 
s^1_{i,i}      \\
s^2_{i,i}      \\
\ldots        \\ 
s^{l-1}_{i,i} 
\end{pmatrix}
-\frac{|f^{-1}(i)|^2}{|\VV{G}|}
\begin{pmatrix}
1    \\
\lambda_0  \\ 
\lambda^2_0 \\
\ldots        \\ 
\lambda^{l-1}_0 
\end{pmatrix}}
\]
has a unique solution $(a_1,\ldots,a_l)$. Moreover, $a_j \cdot |\VV{G}|$ is a non-negative integer for $j=1,\ldots,l$.
\end{theoreman}
\begin{proof}
The matrix of the system is a transposition of Vandermonde matrix, so the determinant is not equal to $0$. Hence the system has a unique solution. Let $f_i=\chi_{f^{-1}(i)}$ be the characteristic function of color $i$. 
By Lemma~\ref{l:eigenvalues} eigenvalues $\lambda_0,\ldots, \lambda_l$ are eigenvalues of $G$.
It is known that $f_i$ belongs to the direct sum of the eigensubspaces corresponding to the eigenvalues of $S$, i.e., $f_i \in U_{\lambda_0} \oplus \ldots \oplus U_{\lambda_l}$ 
(see for example \cite[Property~9]{Tar:perfstruct}). 
So, 
$$f_i=\sum\limits_{z: \varphi_z \in U_{\lambda_{1}}}\alpha_z\varphi_z+
\ldots+ 
\sum\limits_{z: \varphi_z \in U_{\lambda_{l}}}\alpha_z\varphi_z+\alpha_{\overline{0}} \varphi_{\overline{0}},$$
where $\alpha_z$, $z \in \VV{G}$, are complex coefficients. Therefore,

$$M^tf_i=\sum\limits_{z: \varphi_z \in U_{\lambda_{1}}}\lambda^t_1\alpha_z\varphi_z+
\ldots+ 
\sum\limits_{z: \varphi_z \in U_{\lambda_{l}}}\lambda^t_l\alpha_z\varphi_z+\lambda^t_0\alpha_{\overline{0}} \varphi_{\overline{0}}.$$
This representation implies the following relation for $t=0,1,\ldots$: 

$$(M^tf_i,f_i)=\lambda^t_1\sum\limits_{z: \varphi_z \in U_{\lambda_{1}}}|\alpha_z|^2+
\ldots+ 
\lambda^t_{l}\sum\limits_{z: \varphi_z \in U_{\lambda_{l}}}|\alpha_z|^2+
\lambda^{t}_0|\alpha_{\overline{0}}|^2.$$
Since the basis is orthonormal, we have $\alpha_{\overline{0}}=(f_i,\varphi_{\overline{0}})=\frac{|f^{-1}(i)|}{|\VV{G}|^{1/2}}$, 
and hence $|\alpha_{\overline{0}}|^2=\frac{|f^{-1}(i)|^2}{|\VV{G}|}$.
Since $(M^tf_i,f_i)=|f^{-1}(i)| \cdot s^t_{i,i}$ by Lemma~\ref{l:mffs}, it is 
straightforward that $(a_1,\ldots,a_l)$, where $a_j=\sum\limits_{z: \varphi_z \in U_{\lambda_{j}}}|\alpha_z|^2$, is the solution of the system. 
On the other hand, as the basis $\{\varphi_z:z \in \VV{G}\}$ is orthonormal, we have $\alpha_z=(f_i,\varphi_z)$ for any $z \in \VV{G}$. Let us consider subcases.

If $q=2$, then for any $z \in \mathbb Z^n_2$ the function $\varphi_z$ has two distinct values: 
$\frac{\pm 1}{2^{n/2}}$. In this case, $\alpha_z=(f_i,\varphi_z)=\frac{r}{2^{n/2}}$ and $|\alpha_z|^2=\frac{r^2}{2^n}$ for some integer $r$.
Hence $a_j 2^n=\sum\limits_{z: \varphi_z \in U_{\lambda_{j}}}2^n|\alpha_z|^2$ is a non-negative integer.

If $q=3$, then $\varphi_z$ has three distinct values: 
$\frac{-1+\sqrt{3}\sqrt{-1}}{2 \cdot 3^{n/2}},\frac{-1-\sqrt{3}\sqrt{-1}}{2 \cdot 3^{n/2}},\frac{1}{3^{n/2}}$. 
In this case, $\alpha_z=\frac{a+b\sqrt{3}\sqrt{-1}}{2 \cdot 3^{n/2}}$, were $a$ and $b$ are integers and, moreover, they have the same parity. So $|\alpha_z|^2=\frac{a+3b^2}{4\cdot3^{n}}=\frac{r}{3^n}$ for some integer $r$. 
So $a_j 3^n=\sum\limits_{z: \varphi_z \in U_{\lambda_{j}}}3^n|\alpha_z|^2$ is a non-negative integer.

If $G$ is $D(m,n)$, then $\varphi_z$ has four distinct values: 
$\frac{\pm 1}{4^{(2m+n)/2}},\frac{\pm \sqrt{-1}}{4^{(2m+n)/2}}$. 
So $\alpha_z=(f_i,\varphi_z)=\frac{a+b\sqrt{-1}}{4^{(2m+n)/2}}$ for some integers $a$ and $b$. 
Hence $|\alpha_z|^2=\frac{r}{4^{2m+n}}$ for some integer $r$. Hence $a_j 4^{2m+n}=\sum\limits_{z: \varphi_z \in U_{\lambda_{j}}}4^{2m+n}|\alpha_z|^2$ is a non-negative integer. 
\end{proof}

\section{Extended perfect codes are completely regular}\label{s:extarecrg}

\begin{theoreman}\label{t:extiscomp}
\begin{enumerate}
    \item A code $C$ in $H(n,q)$ is extended $1$-perfect if and only if
     $C$ is completely regular with quotient matrix
\[\displaystyle{\begin{pmatrix}
0 \ \ \  & n(q-1) & 0 \\
1 \ \ \  & q-2 & (n-1)(q-1) \\ 
0 \ \ \ & n & n(q-2)
\end{pmatrix}}.\]
    \item A code $C$ in $D(m,n)$ is extended $1$-perfect if and only if
     $C$ is completely regular with quotient matrix
    \[\displaystyle{\begin{pmatrix}
0 \ \ \  & 6m+3n & 0 \\
1 \ \ \  & 2 & 6m+3n-3 \\ 
0 \ \ \ & 2m+n & 4m+2n
\end{pmatrix}}.\]
\end{enumerate}
\end{theoreman}
\begin{proof}
In most parts, the proof for $D(m,n)$ is similar to the proof for $H(2m+n,4)$.  
So we mainly focus on Hamming graphs, and we consider Doob graphs only in the cases for which the proof is different.
Let $C$ be an extended $1$-perfect code in $H(n,q)$ ($D(m,n)$).
Let $f$ be the distance coloring of $H(n,q)$ with respect to $C$, i.e. $f(x)=\min\limits_{y \in C} \{d(x,y)\}$, $x \in \VV{H(n,q)}$.
Since the projection of $C$ in any position is a $1$-perfect code that has the covering radius $1$, 
the covering radius of $C$ equals $2$ hence the set of colors is $\{0,1,2\}$.

Define the following $s^i_j: f^{-1}(i) \to \mathbb Z$,
where $s^i_j(x)$ is the number of vertices of color $j$ in the neighbourhood of $x$, if $f(x)=i$, and otherwise is not defined. 
So, $f$ is a perfect coloring if and only if $s^i_j$ is constant for all $i,j \in \{0,1,2\}$.
Obviously, $s^0_0 \equiv 0$ (as the code distance is $4$), and $s^0_2 \equiv 0$, $s^2_0 \equiv 0$ (by the definition).

Let $y$ be an arbitrary vertex of color $1$. Let us count the values $s^1_0(y)$ and $s^1_1(y)$.
On the one hand, $s^1_0(y) \ge 1$ by the definition. 
On the other hand, $s^1_0(y) \le 1$ (otherwise we have a contradiction with the code distance).
Hence $s^1_0 \equiv 1$. 
Therefore, for any vertex $x$ of color $1$ we can denote by $o(x)$ the unique neighbour of $x$ that has color $0$.
Any vertex $y'$ of color $1$ that adjacent to $y$ belongs to the neighbourhood of $o(y)$ 
(indeed, if $o(y) \ne o(y')$, then $d(o(y),o(y')) \le 3$ that contradicts the code distance).
Therefore, all neighbours of $y$ that have color $1$ belong to the neighbourhood of $o(y)$.
The number of common neighbours of two arbitrary adjacent vertices in a distance-regular graph is uniquely determined by the intersection array.
For $H(n,q)$, it is equal to $q-2$ and for $D(m,n)$, to $2$.
Hence $s^1_0 \equiv 1$ and $s^1_1 \equiv q-2$ ($s^1_1 \equiv 2$ for a Doob graph).
For each vertex $x \in \VV{H(n,q)}$, we have $s^i_0(x)+s^i_1(x)+s^i_2(x)=n(q-1)$, where $i$ is the color of $x$.
Therefore, $s^0_1 \equiv n(q-1)$ and $s^1_2 \equiv (n-1)(q-1)$.

It remains to prove that $s^2_1 \equiv n$ ($s^2_1 \equiv 2m+n$ for $D(m,n)$).
An edge $\{v,u\}$ is called an \emph{$(i,j)$-edge} if $v$ has color $i$ and $u$ has color $j$, or vice versa. 
Denote $\displaystyle{\alpha=\sum\limits_{x \in f^{-1}(2)} s^2_1(x)}$, 
i.e. the number of $(1,2)$-edges.

Let us calculate the values $|f^{-1}(0)|$, $|f^{-1}(1)|$ and $|f^{-1}(2)|$. 
The first value is equal to the cardinality of a $1$-perfect code in $H(n-1,q)$, i.e. $\displaystyle{\frac{q^{n-1}}{(n-1)(q-1)+1}}$.
From the counting of the number of $(0,1)$-edges, we have
$|f^{-1}(1)|=|f^{-1}(0)| n(q-1)=\displaystyle{\frac{n(q-1)q^{n-1}}{(n-1)(q-1)+1}}$. 
Counting the number of $(1,2)$-edges, we find
$\alpha = (n-1)(q-1)|f^{-1}(1)|$.
On the other hand, 
\begin{multline*}
|f^{-1}(2)|= 
q^n-|f^{-1}(0)|-|f^{-1}(1)|=\\
q^{n-1}\frac{q((n-1)(q-1)+1)-n(q-1)-1}{(n-1)(q-1)+1}= \\
q^{n-1}\frac{(n-1)(q-1)^2}{(n-1)(q-1)+1}
\end{multline*}

Hence the average value of $s^2_1$ equals $n$, i.e. 
$\displaystyle{\frac{\alpha}{|f^{-1}(2)|}=n}$ 
(or $2m+n$ for $D(m,n)$).

Let $v$ be a vertex of color $2$ in $H(n,q)$.
The induced subgraph on the set of its neighbours has $n$ connected components, and every component is a $(q-1)$-clique. Hence
$s^2_1(v) \le n$ (otherwise there are two vertices $u$ and $w$ of color $1$ in the same component, 
but all their common neighbours except $v$ also belongs to this component and one of them is $o(u)$, which has color $0$).
Since the average value of $s^2_1$ equals $n$, we have $s^2_1 \equiv n$.

Let $v=(x_1,\ldots,x_m;y_1,\ldots,y_n)$ be a vertex of color $2$ in $D(m,n)$.
Denote by $h_{j,v}$ the induced subgraph on the set 
$\{(x_1,\ldots,x_m;y_1,\ldots,y_{j-1},b,y_{j+1},\ldots,y_n): b \in \mathbb Z_4\}$. This graph is the complete graph $K_4$.
Denote by $d_{i,v}$ the induced subgraph on the vertex set 
$\{(x_1,\ldots,x_{i-1},a,x_{i+1},\ldots,x_m;y_1,\ldots,y_n): a \in \mathbb Z^2_4\}$. This graph is the Shrikhande graph.
Denote by $\alpha_{i,v}$ the number of $(1,2)$-edges in $d_{i,v}$ divided by the number of vertices of color $2$ in $d_{i,v}$. 
Let us prove that for any $i \in \{1,\ldots,m\}$ and $v \in \VV{D(m,n)}$ it follows that $\alpha_{i,v} \le 2$; moreover, if $\alpha_{i,v}=2$, then any vertex of color $2$ in $d_{i,v}$ has exactly two neighbours of color $1$ in $d_{i,v}$.

Let $i \in \{1,\ldots,m\}$ and $v \in \VV{D(m,n)}$. Consider two cases.
If $d_{i,v}$ contains a vertex $u$ of color $0$, 
then $\alpha_{i,v}=2$. 
Indeed, all neighbours of $u$ have color $1$ and other $9$ vertices have color $2$ 
(if some vertex $w$ is at distance $2$ from some vertex of color $0$, 
then $f(w)=2$; otherwise we have a contradiction with the code distance). 
So any vertex of color $2$ has two neighbours of color $1$ (because the Shrikhande graph is strongly regular with parameters $(16,6,2,2)$).
In the second case, there are no vertices of color $0$ in $d_{i,v}$. 
Then the vertices of color $1$ form an independent set 
(indeed, if some vertices $u$ and $w$ are adjacent, 
then $o(u)$ is their common neighbour, but these vertices have only two common neighbours, which also belong to $d_{i,v}$).
So $\alpha_{i,v}=\frac{6x}{16-x}$, where $x$ is the number of vertices of color $1$. 
A maximum independent set in the Shrikhande graph has cardinality $4$; 
moreover, the characteristic function of a maximum independent set is a perfect coloring,
where any vertex that does not belong to this set is adjacent to $2$ vertices from this set (see \cite[Section~2]{BesKro:mdsdoob}).
Hence $\alpha_{i,v} \le 2$; moreover, if $\alpha_{i,v}=2$ ($x=4$), then any vertex of color $2$ has exactly two neighbours of color $1$ in $d_{i,v}$.
As before, for any $j\in \{1,\ldots,n\}$ and $v \in \VV{D(m,n)}$ any vertex of color $2$ has $0$ or $1$ neighbours of color $1$ in the graph $h_{j,v}$. 
Since any $(1,2)$-edge in $D(m,n)$ belongs to exactly one subgraph among the subgraphs $d_{i,v}$ and $h_{j,v}$, 
where $v \in \VV{D(m,n)}$, $i=1,\ldots,m$, $j=1,\ldots,n$, we have 
$\displaystyle{\frac{\alpha}{|f^{-1}(2)|} \le 2m+n}$. 
Moreover, if $\displaystyle{\frac{\alpha}{|f^{-1}(2)|} = 2m+n}$, then $\alpha_{i,v}=2$
for any $i \in \{1,\ldots,m\}$ and $v \in \VV{D(m,n)}$. Hence $s^2_1 \equiv 2m+n$.  

Let us prove the converse statement for the Hamming graphs (for the Doob graphs the proof is similar). 
Let $f$ be a perfect $3$-coloring with quotient matrix $S$ from the theorem statement.
Let us prove that the code $C=f^{-1}(0)$ is an extended $1$-perfect code. 
Since $s_{0,0}=0$, the code distance $d$ is at least $2$. 
Moreover, $s_{1,0}=1$ implies $d \ge 3$.
Suppose that there are different vertices $x$ and $y$ in $C$ such that $d(x,y)=3$. 
In this case, there is a path $(x,v_1,v_2,y)$. 
The vertices $v_1$ and $v_2$ have color $1$. 
Since $v_1$ ($v_2$) has $q-2=s_{1,1}$ common neighbours with $x$ ($y$), we have a contradiction with the fact that $v_1$ and $v_2$ are adjacent.
So $d \ge 4$.
By the counting of $(i,j)$-edges, we have $s_{i,j}|f^{-1}(i)|=s_{j,i}|f^{-1}(j)|$ for any $i, j$. 
It implies $q^n=|C|(1+n(q-1)+(n-1)(q-1)^2)=|C|q((n-1)(q-1)+1)$. 
Therefore, a projection of $C$ in any position is a code with the code distance not less than $3$, whose cardinality achieves the sphere-packing bound. 
So $C$ is an extended $1$-perfect code.
\end{proof}

The following lemma can be checked directly.

\begin{lemman}\label{l:eigen}
The matrix \[\displaystyle{\begin{pmatrix}
0 \ \ \  & n(q-1) & 0 \\
1 \ \ \  & q-2 & (n-1)(q-1) \\ 
0 \ \ \ & n & n(q-2)
\end{pmatrix}}\]
has the following eigenvalues: $\lambda_1=(q-2)$, $\lambda_2=-n$, and $\lambda_0=n(q-1)$.
\end{lemman}

\section{The non-existence of some extended perfect codes}\label{s:nonexist}

Now we can apply Theorem~\ref{t:nessesary} to prove the non-existence of ternary and quaternary extended $1$-perfect codes.

\begin{predln}\label{p:extnonexist}
\begin{enumerate}
\item Let $C$ be an extended $1$-perfect code in $H(n,3)$, where $n=\frac{3^l+1}{2}$, $l \in \mathbb N$. Then $l \le 2$. 
\item Let $C$ be an extended $1$-perfect code in $D(m,n)$ (including the case $D(0,n)=H(n,4)$), where $2m+n=\frac{4^l+2}{3}$, $l \in \mathbb N$. Then $l \le 3$.
\end{enumerate}
\end{predln}
\begin{proof}
1) Let $C$ be an extended $1$-perfect code in $H(n,3)$, where $n=\frac{3^l+1}{2}$ for some positive integer $l$.
The cardinality of $C$ is equal to $3^{n-l-1}$.
By Theorem~\ref{t:extiscomp} and Lemma~\ref{l:eigen} the distance coloring with respect to $C$ is a perfect coloring with quotient matrix,
which has eigenvalues: $\lambda_1=1$, $\lambda_2=-n$ and $\lambda_0=2n$.
Let us consider the system of equations from Theorem~\ref{t:nessesary}

\[\begin{cases}
a_1+a_2=3^{n-l-1}-3^{n-2l-2}\\
a_1-na_2=-2n \cdot 3^{n-2l-2}.
\end{cases}\]
From this system we have

$$\displaystyle{a_2 \cdot 3^n=\frac{3^{2n-2l-2}(3^{l+1}+3^l)}{\frac{3^l+3}{2}}=\frac{3^{2n-l-3}2^3}{3^{l-1}+1}}.$$

By Theorem~\ref{t:nessesary} the number $a_2 \cdot 3^n$ is integer.
Since the denominator $3^{l-1}+1$ and $3^{2n-l-3}$ are relatively prime, it follows that $3^{l-1}+1$ is a divisor of $8$. 
This implies $l=1$ or $l=2$.

2) Let $C$ be an extended $1$-perfect code in $D(m,n)$, where $2m+n=\frac{4^l+2}{3}$ for some positive integer $l$.

In this case, we have the following system of equations
\[\begin{cases}
a_1+a_2=4^{2m+n-l-1}-4^{2m+n-2l-2}\\
2a_1-(2m+n)a_2=-3(2m+n)4^{2m+n-2l-2}.
\end{cases}\]

From this system we have

$$\displaystyle{a_2 \cdot 4^{2m+n}=\frac{4^{4m+2n-2l-2}(2 \cdot 4^{l+1}-2+4^l+2)}{\frac{4^l+8}{3}}=\frac{4^{4m+2n-l-3}3^3}{4^{l-1}+2}}.$$

By Theorem~\ref{t:nessesary} the number $a_2 \cdot 4^{2m+n}$ is integer. 
If $l=1$, then $4^2 \cdot a_2=9$. Let $l>1$. 
Since the greatest common divisor of the denominator $4^{l-1}+2$ and $4^{4m+2n-l-3}$ equals $2$, 
it follows that $2 \cdot 4^{l-2}+1$ divides $27$.
This implies $l \in \{2,3\}$. So $l \le 3$.
\end{proof}

The two following propositions solve the remaining cases in $H(n,3)$ and $D(m,n)$, 
and codes of odd length in $H(n,q)$ for all $q$.
The proofs of these propositions are  particular cases of the method described in \cite{Kro:struct}.

\begin{predln}\label{p:odddiam}
Let $C$ be an extended $1$-perfect code in $H(n,q)$. Then $n$ is even.
\end{predln}
\begin{proof}
Let $C$ be an extended $1$-perfect code in $H(n,q)$ and $f$ be the distance coloring with respect to $C$. 
Consider an arbitrary vertex $a$ of color $2$. 
Denote by $W^i_j$ the set of vertices of color $i$ that are at the distance $j$ from $a$ 
and denote $W_j=W^0_j \cup W^1_j \cup W^2_j$.  
On the one hand, any vertex $x \in W^1_1$ is adjacent to exactly $1$ vertex from $W^0_2$. 
On the other hand, any vertex $y \in W^0_2$ has $2$ neighbours in $W_1$ and they have color $1$.
Hence $|W^1_1|=2|W^0_2|$, and so $|W^1_1|$ is even.
But from Theorem~\ref{t:extiscomp} we have $|W^1_1|=n$.
\end{proof}

Recall that a code $C$ in $H(n,q)$ is called an \emph{MDS code} with distance $d$ if its cardinality 
achieves the Singleton bound, i.e. $|C|=q^{n-d+1}$. 
In the case $n=q+2$, the definitions of an extended $1$-perfect code and an MDS code with distance $4$ are equivalent.

\begin{coroll}
If $q$ is odd, then there are no MDS codes with distance $4$ in $H(q+2,q)$. 
\end{coroll}

\begin{coroll}
Let $q=p^m$ be an odd prime power, and let $C$ be an extended $1$-perfect code in $H(n,q)$.
Then $\displaystyle{n=\frac{q^{l}+q-2}{q-1}}$ for some odd $l$. 
\end{coroll}

\begin{predln}\label{p:remaincases}
There are no extended $1$-perfect codes in $D(m,n)$, where $2m+n=22$.
\end{predln}
\begin{proof}
Let $C$ be an extended $1$-perfect code in $D(m,n)$, where $2m+n=22$, and let $f$ be the distance coloring with respect to $C$. 
Consider an arbitrary vertex $a$ of color $2$.
Denote by $W^i_j$ the set of vertices of color $i$ that are at the distance $j$ from $a$ 
and denote $W_j=W^0_j \cup W^1_j \cup W^2_j$.  
By Theorem~\ref{t:extiscomp} we have $|W^0_1|=0$, $|W^1_1|=22$, and $|W^2_1|=44$. 
As in proof of Proposition~\ref{p:odddiam}, we have $2|W^0_2|=|W^1_1|$, so $|W^0_2|=11$. 
Let us count the number $w$ of edges $(x,y)$ such that $x \in W_1$ and $y \in W^1_2$. 
This number is equal to $(22 \cdot 2 + 44 \cdot 22 - 2t-r)$, 
where $t$ is the number of $(1,1)$-edges and $r$ is the number of $(1,2)$-edges 
in the induced subgraph on the set of vertices $W_1$. 
It follows from the intersection array that this subgraph is $2$-regular, and hence $2t+r=2|W^1_1|=44$.
So $w=22 \cdot 2+44 \cdot 22 - 44=968$. 
On the other hand, $w=2|W^1_2|$, so $|W^1_2|=484$.
Let us count the number of $(0,1)$-edges that are incident to some vertex from $W^1_2$.
This number is equal to $|W^1_2|=484$; on the other hand, it is equal to $6|W^0_2|+3|W^0_3|=66+3|W^0_3|$.
We find that $3|W^0_3|=418$. Since $|W^0_3|$ is integer, we have a contradiction.
\end{proof}

Remind that formally the singleton from any vertex in $H(2,3)$,  $D(0,2)$ or $D(1,0)$ is an extended $1$-perfect code, called trivial.
Also all extended $1$-perfect codes in $D(m,n)$, where $2m+n=6$, are characterized in \cite{Alderson:MDS4, BesKro:mdsdoob}. 
From Propositions~\ref{p:extnonexist}, \ref{p:odddiam}, and \ref{p:remaincases} we have the following statement.

\begin{theoreman}\label{t:parameters}
\begin{enumerate}
\item An extended $1$-perfect code in $H(n,3)$ exists if and only if $n=2$.
\item An extended $1$-perfect code in $D(m,n)$ (including the case $D(0,n)=H(n,4)$) exists if and only if 
$(m,n)=(0,2)$, or $(m,n)=(1,0)$, or $(m,n)=(0,6)$, or $(m,n)=(2,2)$.
\item For any $q$, there are no extended $1$-perfect codes in $H(n,q)$ if $n$ is odd.
\end{enumerate}
\end{theoreman}

\section*{Acknowledgements}

The author is grateful to Denis Krotov, Vladimir Potapov, and Ev Sotnikova for helpful remarks and introducing him to some background.

\bibliographystyle{unsrt}
\bibliography{k}

\providecommand\href[2]{#2} \providecommand\url[1]{\href{#1}{#1}}
  \def\DOI#1{{\small {DOI}:
  \href{http://dx.doi.org/#1}{#1}}}\def\DOIURL#1#2{{\small{DOI}:
  \href{http://dx.doi.org/#2}{#1}}}
\begin{thebibliography}{10}

\bibitem{Tar:perfstruct}
A.~A. Taranenko.
\newblock Algebraic properties of perfect structures.
\newblock E-print 1906.10430v2, arXiv.org, 2020.
\newblock Available at \url{https://arxiv.org/abs/1906.10430v2}.

\bibitem{FDF:CorrImmBound}
D.~G. Fon-Der-Flaass.
\newblock A bound on correlation immunity.
\newblock {\em \href{http://semr.math.nsc.ru}{Sib. Ehlektron. Mat. Izv.}},
  4:133--135, 2007.
\newblock Online: \url{http://mi.mathnet.ru/eng/semr149}.

\bibitem{Bierbrauer:95}
J.~Bierbrauer.
\newblock Bounds on orthogonal arrays and resilient functions.
\newblock {\em
  \href{http://onlinelibrary.wiley.com/journal/10.1002/(ISSN)1520-6610}{J.
  Comb. Des.}}, 3(3):179--183, 1995.
\newblock \DOI{10.1002/jcd.3180030304}.

\bibitem{Friedman:92}
J.~Friedman.
\newblock On the bit extraction problem.
\newblock In {\em Foundations of Computer Science, IEEE Annual Symposium on},
  pages 314--319, Los Alamitos, CA, USA, 1992. IEEE Computer Society.
\newblock \DOI{10.1109/SFCS.1992.267760}.

\bibitem{Pot:2012:color}
V.~N. Potapov.
\newblock On perfect $2$-colorings of the $q$-ary $n$-cube.
\newblock {\em
  \href{http://www.sciencedirect.com/science/journal/0012365X}{Discrete
  Math.}}, 312(6):1269--1272, 2012.
\newblock \DOI{10.1016/j.disc.2011.12.004}.

\bibitem{Pot:2010:correng}
V.~N. Potapov.
\newblock On perfect colorings of boolean $n$-cube and correlation immune
  functions with small density.
\newblock {\em \href{http://semr.math.nsc.ru}{Sib. Ehlektron. Mat. Izv.}},
  7:372--382, 2010.
\newblock In Russian, English abstract.

\bibitem{Kro:OA1536}
D.~S. Krotov.
\newblock On the {O}{A}(1536,13,2,7) and related orthogonal arrays.
\newblock {\em Discrete Mathematics}, 343(2):111659, 2020.
\newblock \DOI{10.1016/j.disc.2019.111659}.

\bibitem{Kro:2m-3}
D.~S. Krotov.
\newblock On the binary codes with parameters of doubly-shortened $1$-perfect
  codes.
\newblock {\em \href{http://link.springer.com/journal/10623}{Des. Codes
  Cryptography}}, 57(2):181--194, 2010.
\newblock \DOI{10.1007/s10623-009-9360-5}.

\bibitem{Kro:2m-4}
D.~S. Krotov.
\newblock On the binary codes with parameters of triply-shortened $1$-perfect
  codes.
\newblock {\em \href{http://link.springer.com/journal/10623}{Des. Codes
  Cryptography}}, 64(3):275--283, 2012.
\newblock \DOI{10.1007/s10623-011-9574-1}.

\bibitem{MWS}
F.~J. MacWilliams and N.~J.~A. Sloane.
\newblock {\em The Theory of Error-Correcting Codes}.
\newblock Amsterdam, Netherlands: North Holland, 1977.

\bibitem{GolombPosner64}
S.~W. Golomb and E.~C. Posner.
\newblock Rook domains, latin squares, and error-distributing codes.
\newblock {\em
  \href{http://ieeexplore.ieee.org/xpl/RecentIssue.jsp?punumber=18}{IEEE Trans.
  Inf. Theory}}, 10(3):196--208, 1964.
\newblock \DOI{10.1109/TIT.1964.1053680}.

\bibitem{HedRoos:2011}
O.~Heden and C.~Roos.
\newblock The non-existence of some perfect codes over non-prime power
  alphabets.
\newblock {\em
  \href{http://www.sciencedirect.com/science/journal/0012365X}{Discrete
  Math.}}, 311(14):1344--1348, 2011.
\newblock \DOI{10.1016/j.disc.2011.03.024}.

\bibitem{Lenstra72}
H.~W. Lenstra, Jr.
\newblock Two theorems on perfect codes.
\newblock {\em
  \href{http://www.sciencedirect.com/science/journal/0012365X}{Discrete
  Math.}}, 3(1-3):125--132, 1972.
\newblock \DOI{10.1016/0012-365X(72)90028-3}.

\bibitem{Lloyd}
S.~P. Lloyd.
\newblock Binary block coding.
\newblock {\em \href{http://www.alcatel-lucent.com/bstj/}{Bell Syst. Tech.
  J.}}, 36(2):517--535, 1957.
\newblock \DOI{10.1002/j.1538-7305.1957.tb02410.x}.

\bibitem{Kro:pfdoob}
D.~S. Krotov.
\newblock The existence of perfect codes in {D}oob graphs.
\newblock {\em IEEE Transactions on Information Theory}, 66(3):1423--1427,
  2020.
\newblock \DOI{10.1109/TIT.2019.2946612}.

\bibitem{Kro:perfect-doob}
D.~S. Krotov.
\newblock Perfect codes in {D}oob graphs.
\newblock {\em \href{http://link.springer.com/journal/10623}{Des. Codes
  Cryptography}}, 80(1):91--102, 2016.
\newblock \DOI{10.1007/s10623-015-0066-6}.

\bibitem{SHK:addperfdoob}
M.~Shi, D.~Huang, and D.~Krotov.
\newblock Additive perfect codes in {D}oob graphs.
\newblock {\em \href{http://link.springer.com/journal/10623}{Des. Codes
  Cryptography}}, 87(8):1857--1869, 2019.
\newblock \DOI{10.1007/s10623-018-0586-y}.

\bibitem{Delsarte:1973}
P.~Delsarte.
\newblock {\em An Algebraic Approach to Association Schemes of Coding Theory},
  volume~10 of {\em Philips Res. Rep., Supplement}.
\newblock N.V. Philips' Gloeilampenfabrieken, Eindhoven, Netherlands, 1973.

\bibitem{Neu:crg}
A.~Neumaier.
\newblock Completely regular codes.
\newblock {\em Discrete Mathematics}, 106-107:353--360, 1992.
\newblock \DOI{10.1016/0012-365X(92)90565-W}.

\bibitem{BRZ:crg}
J.~Borges, J.~Rif\`a, and V.~A. Zinoviev.
\newblock On completely regular codes.
\newblock {\em Problems of Information Transmission}, 55(1):1--45, 2019.
\newblock \DOI{10.1134/S0134347519010017}.

\bibitem{FDF:PerfCol}
D.~G. Fon-Der-Flaass.
\newblock Perfect $2$-colorings of a hypercube.
\newblock {\em \href{http://link.springer.com/journal/11202}{Sib. Math. J.}},
  48(4):740--745, 2007.
\newblock \DOI{10.1007/s11202-007-0075-4} translated from
  \href{http://www.mathnet.ru/php/journal.phtml?jrnid=smj\&option_lang=eng}{Sib.
  Mat. Zh.} 48(4) (2007), 923-930.

\bibitem{BKMTV:perfcolinham}
E.~Bespalov, D.~Krotov, A.~Matiushev, A.~Taranenko, and K.~Vorob{'}ev.
\newblock Perfect $2$-colorings of {H}amming graphs.
\newblock E-print 1911.13151v2, arXiv.org, 2020.
\newblock Available at \url{https://arxiv.org/abs/1911.13151v2}.

\bibitem{KKM:smallval}
J.~Koolen, D.~Krotov, and W.~Martin.
\newblock Completely regular codes: tables.
\newblock \url{https://sites.google.com/site/completelyregularcodes/}.

\bibitem{AhlAydKha}
R.~Ahlswede, H.~K. Aydinian, and L.~H. Khachatrian.
\newblock On perfect codes and related concepts.
\newblock {\em \href{http://link.springer.com/journal/10623}{Des. Codes
  Cryptography}}, 22(3):221--237, 2001.
\newblock \DOI{10.1023/A:1008394205999}.

\bibitem{Hill:Caps}
R.~Hill.
\newblock Caps and codes.
\newblock {\em
  \href{http://www.sciencedirect.com/science/journal/0012365X}{Discrete
  Math.}}, 22:111--137, 1978.
\newblock \DOI{10.1016/0012-365X(78)90120-6}.

\bibitem{Alderson:MDS4}
T.~L. Alderson.
\newblock $(6,3)$-{MDS} codes over an alphabet of size $4$.
\newblock {\em \href{http://link.springer.com/journal/10623}{Des. Codes
  Cryptography}}, 38(1):11--40, 2006.
\newblock \DOI{10.1007/s10623-004-5659-4}.

\bibitem{BesKro:mdsdoob}
E.~Bespalov and D.~Krotov.
\newblock {MDS} codes in {D}oob graphs.
\newblock {\em Problems of Information Transmission}, 53:136--154, 2017.
\newblock \DOI{10.1134/S003294601702003X}.

\bibitem{Ball:2012:1}
S.~Ball.
\newblock On sets of vectors of a finite vector space in which every subset of
  basis size is a basis.
\newblock {\em \href{http://www.ems-ph.org/journals/journal.php?jrn=jems}{J.
  Eur. Math. Soc.}}, 14(3):733--748, 2012.
\newblock \DOI{10.4171/JEMS/316}.

\bibitem{KKO:smallMDS}
J.~I. Kokkala, D.~S. Krotov, and P.~R.~J. {\"O}sterg{\aa}rd.
\newblock On the classification of {MDS} codes.
\newblock {\em
  \href{http://ieeexplore.ieee.org/xpl/RecentIssue.jsp?punumber=18}{IEEE Trans.
  Inf. Theory}}, 61(12):6485--6492, December 2015.
\newblock \DOI{10.1109/TIT.2015.2488659}.

\bibitem{KokOst:further}
J.~I. Kokkala and P.~R.~J. {\"O}sterg{\aa}rd.
\newblock Further results on the classification of {MDS} codes.
\newblock {\em \href{http://aimsciences.org/journals/amc/}{Adv. Math.
  Commun.}}, 10(3):489--498, August 2016.
\newblock \DOI{10.3934/amc.2016020}.

\bibitem{GST:newupbounds}
D.~Gijswijt, A.~Schrijver, and H.~Tanaka.
\newblock New upper bounds for nonbinary codes based on the {T}erwilliger
  algebra and semidefinite programming.
\newblock {\em Journal of Combinatorial Theory, Series A}, 113(8):1719--1731,
  2006.
\newblock \DOI{10.1016/j.jcta.2006.03.010}.

\bibitem{Godsil93}
C.~D. Godsil.
\newblock {\em Algebraic Combinatorics}.
\newblock Chapman and Hall, New York, 1993.

\bibitem{Kro:struct}
D.~S. Krotov.
\newblock On weight distributions of perfect colorings and completely regular
  codes.
\newblock {\em \href{http://link.springer.com/journal/10623}{Des. Codes
  Cryptography}}, 61(3):315--329, 2011.
\newblock \DOI{10.1007/s10623-010-9479-4}.

\end{thebibliography}

\end{document}